\documentclass[10pt]{article}
\usepackage{amsmath,amssymb,amsthm,amscd}
\numberwithin{equation}{section}

\newtheorem{prop}{Proposition}[section]
\newtheorem{theo}[prop]{Theorem}

\newtheorem{rema}[prop]{Remark}

\def\begeq{\begin{equation}}
\def\endeq{\end{equation}}

\def\<{\langle}
\def\>{\rangle}

\begin{document}

\title{Scalar Curvature Bound for K\"ahler-Ricci Flows over 
Minimal Manifolds of General Type}

\author{Zhou Zhang \\
University of Michigan, at Ann Arbor} 
\date{}
\maketitle

\section{Introduction and setup}

We consider the following K\"ahler-Ricci flow over a closed manifold $X$ 
of complex dimension $n\geqslant 2$,
$$\frac{\partial\tilde\omega_t}{\partial t}=-{\rm Ric}(\tilde\omega_t)-\tilde
\omega_t, ~~~~ \tilde\omega_0=\omega_0$$
where $\omega_0$ is any K\"ahler metric. 

In this short note, we are going to prove the following theorem. Some 
classic computations used during the process might be of more interest. 
Most of them are more or less quoted directly from \cite{song-tian}. 

\begin{theo}
Suppose $X$ is a minimal model of general type, then the K\"ahler-Ricci 
flow above has bounded scalar curvature. 
\end{theo}

As in \cite{t-znote}, define $\omega_t:=\omega_\infty+e^{-t}(\omega_0-
\omega_\infty)$ as the background form for the flowing metric, with 
$\omega_\infty=-{\rm Ric}(\Omega)$ for a smooth volume form over $X$ 
and one has $[\omega_\infty]=-c_1(X)=K_X$ cohomologically. Then 
$[\tilde\omega_t]=[\omega_t]$ and we can assume $\tilde\omega_t=
\omega_t+\sqrt{-1}\partial\bar\partial u$. The following evolution of 
the space-time function $u$ (usually called the metric potential of the 
flowing metric),  
$$\frac{\partial u}{\partial t}={\rm log}\frac{(\omega_t+\sqrt{-1}\partial\bar
\partial u)^n}{\Omega}-u, ~~~~ u(0, \cdot)=0,$$
would imply the metric flow above. They are indeed equivalent to each 
other by the basics on the existence and uniqueness of these flows (as 
summarized in \cite{thesis}, for example). 

Without further clarification, all the constants 
appearing later are positive. The same letter 
might stand for different (but fixed) constants 
at different places.

By direct Maximum Principle argument, this equation above already gives 
$$u\leqslant C$$ 
as long as the flow exists.  

In the following, we summarize some useful computation and estimates 
already known (as in \cite{t-znote} and \cite{thesis}) without any assumption 
on the closed manifold, $X$.

In this note, the Laplacian $\Delta$ and norm $|\cdot|$, below are always with 
respect to the metric along the flow, $\tilde\omega_t$. 
$$\frac{\partial}{\partial t}(\frac{\partial u}{\partial t})=\Delta(\frac{\partial u}
{\partial t})-e^{-t}\<\tilde\omega_t, \omega_0-\omega_\infty\>-\frac{\partial u}
{\partial t},$$
which is just the $t$-derivative of the orginal equation and 
has the following transformation,
$$\frac{\partial}{\partial t}(e^t\frac{\partial u}{\partial t})=\Delta(e^t\frac{\partial u}
{\partial t})-\<\tilde\omega_t, \omega_0-\omega_\infty\>,$$ 
$$\frac{\partial}{\partial t}(\frac{\partial u}{\partial t}+u)=\Delta(\frac{\partial u}
{\partial t}+u)-n+\<\tilde\omega_t, \omega_\infty\>.$$

A proper linear combination of these equations provides 
the following "finite time version" of the previous equation,
$$\frac{\partial}{\partial t}((1-e^{t-T})\frac{\partial u}{\partial t}+u)=\Delta((1-
e^{t-T})\frac{\partial u}{\partial t}+u)-n+\<\tilde\omega_t, \omega_T\>.$$
As $T=\infty$, this naturally gives back the equation above. 

The difference of the $t$-derivative equation and its transformation gives
$$\frac{\partial}{\partial t}((1-e^t)\frac{\partial u}{\partial t}+u)=\Delta((1-e^t)
\frac{\partial u}{\partial t}+u)-n+\<\tilde\omega_t, \omega_0\>,$$
which implies the "essential decreasing" of metric potential along 
the flow, i. e.,  
$$\frac{\partial u}{\partial t}\leqslant\frac{nt+C}{e^t-1}.$$
Notice that this estimate only depends on the initial value 
of $u$ and its upper bound along the flow. It is uniform 
away from the initial time.

Another $t$-derivative gives 
$$\frac{\partial}{\partial t}(\frac{\partial^2 u}{\partial t^2})=\Delta(\frac{\partial^2 u}
{\partial t^2})+e^{-t}\<\tilde\omega_t, \omega_0-\omega_\infty\>-\frac{\partial^2 u}
{\partial t^2}-|\frac{\partial\tilde\omega_t}{\partial t}|^2_{\tilde\omega_t}.$$
Take summation with the first $t$-derivative to arrive at 
$$\frac{\partial}{\partial t}(\frac{\partial^2 u}{\partial t^2}+\frac{\partial u}{\partial t})=
\Delta(\frac{\partial^2 u}{\partial t^2}+\frac{\partial u}{\partial t})-(\frac{\partial^2 u}
{\partial t^2}+\frac{\partial u}{\partial t})-|\frac{\partial\tilde\omega_t}{\partial t}|^2_
{\tilde\omega_t},$$
which gives
$$\frac{\partial^2 u}{\partial t^2}+\frac{\partial u}{\partial t}\leqslant Ce^{-t},$$
which implies the "essential decreasing" of volume form along the flow, i. e. ,   
$$\frac{\partial}{\partial t}(\frac{\partial u}{\partial t}+u)\leqslant Ce^{-t},$$
which also induces
$$\frac{\partial u}{\partial t}\leqslant Ce^{-t}.$$

Let's point out that most of these estimates depends on the initial values of 
$\frac{\partial^2 u}{\partial t^2}$, $\frac{\partial u}{\partial t}$ and $u$, which 
is not a problem here. \\ 

Also rewrite the metric flow equation as follows,
$${\rm Ric}(\tilde\omega_t)=-\sqrt{-1}\partial\bar\partial(u+\frac{\partial u}{\partial t})
-\omega_\infty.$$

Taking trace with respect to $\tilde\omega_t$ for the original metric flow equation 
and the one above, we have
$$R=e^{-t}\<\tilde\omega_t, \omega_0-\omega_\infty\>-\Delta(\frac{\partial u}
{\partial t})-n=-\Delta(u+\frac{\partial u}{\partial t})-\<\tilde\omega_t, \omega_
\infty\>,$$
where $R$ denotes the scalar curvature of $\tilde\omega_t$. Using the equations 
above, we also have
$$R=-n-\frac{\partial}{\partial t}(\frac{\partial u}{\partial t}+u),$$
and so the estimate got for $\frac{\partial}{\partial t}(\frac
{\partial u}{\partial t}+u)$ before is equivalent to the well 
known fact for scalar curvature.

\section{Current interest} 

We have proved in \cite{t-znote} that the flow exists (smoothly) as long as the 
cohomology class $[\tilde\omega_t]=[\omega_t]$ remains to be K\"ahler, which 
can be grasped by simple algebraic concern. As usual, let's define 
$$T=sup\{t|~[\omega_t]~is~Kaehler.\}.$$   

Of course, our main interest is on the case when $[\omega_T]$ is not K\"ahler. 
$T$ can be either infinite or finite here. From now on, $T\leqslant \infty$ unless 
explicitly stated otherwise. We only consider smooth solution of K\"ahler-Ricci 
flow in $[0, T)\times X$. 

At this moment, we focus on the case when the smooth limiting background 
form $\omega_T\geqslant 0$. It is essentially equivalent to assume $[\omega
_T]$ has a smooth non-negative representative (or semi-ample for an algebraic 
geometry background). \\

For $T<\infty$, we have 
$$\frac{\partial}{\partial t}((1-e^{t-T})\frac{\partial u}{\partial t}+u)=\Delta((1-
e^{t-T})\frac{\partial u}{\partial t}+u)-n+\<\tilde\omega_t, \omega_T\>$$
with the "$T$" in the equation chosen to the "$T$" above. With $\omega_T
\geqslant 0$, by Maximum Principle, one has
$$(1-e^{t-T})\frac{\partial u}{\partial t}+u\geqslant -C.$$ 

As $u\leqslant C$ and $\frac{\partial u }{\partial t}\leqslant C$, we can conclude 
that 
$$u\geqslant -C, ~~~~ \frac{\partial u}{\partial t}\geqslant 
-\frac{C}{1-e^{t-T}}\thicksim \frac{C}{t-T}.$$ 

The situation for $T=\infty$ is different. Further assuming $[\omega_\infty]=K_X$ 
is also big, we have the lower bound of $u$, which can actually imply the lower 
bound of $\frac{\partial u}{\partial t}$ (as in \cite{t-znote}, \cite{ey-gu-ze}, \cite{zh} 
and \cite{thesis} by generalizing Kolodziej's $L^\infty$ estimate for complex 
Monge-Amp\`ere equation as summarized in \cite{kojnotes}). More precisely, 
one has 
$$u\geqslant -C, ~~~~ \frac{\partial u}{\partial t}\geqslant 
-C.$$ 

These are the cases under consideration right now. We know that it is always true 
that
$$|(1-e^{t-T})\frac{\partial u}{\partial t}+u|\leqslant C,$$
where the meaning for the case of $T$ being infinity is also natural.

\section{Further computation}

\subsection{Parabolic Schwarz estimate}

Use the following setup as in \cite{song-tian}. Suppose we have a holomorphic 
(non-trivial) map $F:~X \to Y$ between smooth closed complex manifolds, 
$\omega$ is a K\"ahler metric over $Y$ and $\tilde\omega_t$ is metric under 
K\"ahler-Ricci over $X$. Let $\phi=\<\tilde\omega_t, F^*\omega\>$, then one has, 
over $[0, T)\times X$, 
$$(\frac{\partial}{\partial t}-\Delta)\phi\leqslant \phi+C\phi^2-H,$$
where $\Delta$ is with respect to the metric $\tilde\omega_t$, $C$ is related to 
the bisectional curvature bound of $\omega$ over $Y$ and $H\geqslant 0$ is 
described as follows. Using normal coordinates locally over $X$ and $Y$, with 
indices $i, j$ and $\alpha, \beta$, $\phi=|F_i^\alpha|^2$ and $H=|F_{ij}^\alpha|
^2$ with all the summations. Notice that the normal coordinate over $X$ is 
changing along the flow. Furthurmore, one has 
$$(\frac{\partial}{\partial t}-\Delta){\rm log}\phi\leqslant C\phi+1.$$

For our case, the map $F$ is coming from the class $[\omega_T]$ with $Y$ being 
some projective space $\mathbb{CP}^N$, and so $\omega_T$ is $F^*\omega$ 
where $\omega$ is (some mutiple of) Fubini-Study metric over $Y$. \\

Set $v=(1-e^{t-T})\frac{\partial u}{\partial t}+u$ and $|v|\leqslant C$ for our concern from 
before. We have 
$$(\frac{\partial}{\partial t}-\Delta)v=-n+\<\tilde\omega_t, \omega_T\>=-n+\phi.$$

By taking a large enough positive constant $A$, the following inequality is true,
$$(\frac{\partial}{\partial t}-\Delta)({\rm log}\phi-Av)\leqslant -C\phi+C.$$

Since $v$ is bounded, Maximum Principle can be used to deduce $\phi\leqslant C$, 
i.e., 
$$\<\tilde\omega_t, \omega_T\>\leqslant C.$$

As in \cite{song-tian}, in complex surface case, when the corresponding map $F$ 
gives a fiber bundle structure of $X$ with base and fiber spaces being of complex 
dimension $1$, then (at least for the regular part), one has, when restricted to each 
fiber, 
$$\frac{\tilde\omega_t}{\omega_0}=\frac{\tilde\omega_t\wedge\omega_T}{\omega_0
\wedge\omega_T}=\frac{\tilde\omega_t\wedge\omega_T}{\tilde\omega_t^2}\cdot\frac
{\tilde\omega_t^2}{\omega_0\wedge\omega_T}=\frac{1}{2}\<\tilde\omega_t, 
\omega_T\>\cdot\frac{\tilde\omega_t^2}{\omega_0\wedge\omega_T}\leqslant C.$$

These two estimates above give us the picture that $\tilde\omega_t$ will not collapse 
horizontally and will be also bounded fiberwisely as $t\to T$. 

\subsection{Gradient and Laplacian estimates}

In this part, we consider gradient and Laplacian estimates for $v$. Recall that
$$(\frac{\partial}{\partial t}-\Delta)v=-n+\phi, ~~~ \phi=\<\tilde\omega_t, 
\omega_T\>.$$

Standard computation (as in \cite{song-tian}) then gives:
$$(\frac{\partial}{\partial t}-\Delta)(|\nabla v|^2)=|\nabla v|^2-|\nabla\nabla v|^2-
|\nabla\bar\nabla v|^2+2{\rm Re}(\nabla\phi, \nabla v),$$
$$(\frac{\partial}{\partial t}-\Delta)(\Delta v)=\Delta v+({\rm Ric}(\tilde\omega_t), 
\sqrt{-1}\partial\bar\partial v)+\Delta\phi.$$

Again, all the $\nabla$, $\Delta$ and $(\cdot, \cdot)$ are with respect to 
$\tilde\omega_t$ and $\nabla\bar\nabla v$ is just $\partial\bar\partial v$. \\

Consider the quantity $\Psi=\frac{|\nabla v|^2}{C-v}$. Since $v$ 
is bounded, one can easily make sure the denominator 
is positive, bounded and also away from $0$. We have 
the following computation,
\begin{equation}
\begin{split}
&~~ (\frac{\partial}{\partial t}-\Delta)\Psi\\
&= (\frac{\partial}{\partial t}-\Delta)(\frac{|\nabla v|^2}{C-v})\\
&= \frac{1}{C-v}\cdot\frac{\partial}{\partial t}(|\nabla v|^2)+\frac{|\nabla v|^2}{(C-v)^2}\cdot
\frac{\partial v}{\partial t}-\bigl( \frac{(|\nabla v|^2)_{\bar i}}{C-v}+\frac
{v_{\bar i}|\nabla v|^2}{(C-v)^2}\bigr)_i\\
&= \frac{|\nabla v|^2}{(C-v)^2}\cdot(\frac{\partial}{\partial t}-\Delta)v+\frac{1}{C-v}\cdot
(\frac{\partial}{\partial t}-\Delta)(|\nabla v|^2)-\frac{v_i\cdot(|\nabla v|^2)_{\bar i}}{(C-v)^2}-
v_{\bar i}\cdot\bigl(\frac{|\nabla v|^2}{(C-v)^2}\bigr)_i\\
&= \frac{|\nabla v|^2}{(C-v)^2}\cdot(\frac{\partial}{\partial t}-\Delta)v+\frac{1}{C-v}\cdot
(\frac{\partial}{\partial t}-\Delta)(|\nabla v|^2)-\frac{2{\rm Re}(\nabla v, \nabla |\nabla v|^2)}
{(C-v)^2}-\frac{2|\nabla v|^4}{(C-v)^3}. \nonumber
\end{split}
\end{equation} 

Plug in the results from before and rewrite the differential 
equality for $\Psi$,
\begin{equation}
\begin{split}
& ~~~~(\frac{\partial}{\partial t}-\Delta)\Psi\\
&= \frac{(-n+\phi)|\nabla v|^2}{(C-v)^2}+\frac{|\nabla v|^2-|\nabla\nabla v|^2-
|\nabla \bar\nabla v|^2}{C-v}+\frac{2{\rm Re}(\nabla\phi, \nabla v)}{C-v}\\
& ~~~~-\frac{2{\rm Re}(\nabla v, \nabla |\nabla v|^2)}{(C-v)^2}-\frac
{2|\nabla v|^4}{(C-v)^3}.
\end{split}
\end{equation}

The computations below are useful when coming to 
transform the expression above.  
\begin{equation}
\begin{split}
|(\nabla v, \nabla |\nabla v|^2)|
&= |v_i(v_j v_{\bar j})_{\bar i}|\\
&= |v_iv_{\bar j}v_{j\bar i}+v_iv_jv_{\bar j\bar i }|\\
&\leqslant |\nabla v|^2(|\nabla\nabla v|+|\nabla\bar\nabla v|)\\
&\leqslant \sqrt2 |\nabla v|^2(|\nabla\nabla v|^2+|\nabla\bar\nabla v|^2)^
{\frac{1}{2} }. \nonumber
\end{split}
\end{equation}
$$\nabla\Psi=\nabla\bigl(\frac{|\nabla v|^2}{C-v}\bigr)=\frac
{\nabla(|\nabla v|^2)}{C-v}+\frac{|\nabla v|^2\nabla v}{(C-v)^2}.$$

Using also the bounds for $\phi$ and $C-v$, we can have 
the following computation with $\epsilon$ representing small 
positive constant (different from place to place),  
\begin{equation}
\begin{split}
&~~ (\frac{\partial}{\partial t}-\Delta)\Psi\\
&\leqslant C|\nabla v|^2+\epsilon\cdot|\nabla \phi|^2-C(|\nabla\nabla v|^2+|\nabla\bar\nabla v|^2)\\
& ~~~~ -(2-\epsilon){\rm Re}(\nabla \Psi, \frac{\nabla v}{C-v})-\epsilon\cdot\frac{{\rm 
Re}(\nabla v, \nabla |\nabla v|^2)}{(C-v)^2}-\epsilon\cdot\frac{|\nabla v|^4}{(C-v)^3}\\
&\leqslant C|\nabla v|^2+\epsilon\cdot|\nabla \phi|^2-C(|\nabla\nabla v|^2+|\nabla\bar\nabla v|^2)\\
& ~~~~ -(2-\epsilon){\rm Re}(\nabla \Psi, \frac{\nabla v}{C-v})+\epsilon\cdot(|\nabla\nabla v|^2+
|\nabla\bar\nabla v|^2)-\epsilon\cdot|\nabla v|^4\\
&\leqslant C|\nabla v|^2+\epsilon\cdot|\nabla \phi|^2-(2-\epsilon){\rm Re}(\nabla \Psi, \frac
{\nabla v}{C-v})-\epsilon\cdot|\nabla v|^4.\nonumber
\end{split}
\end{equation}

We need a few more calculations to set up Maximum Principle argument. Recall that 
$\phi=\<\tilde\omega_t, \omega_T\>$ and, 
$$(\frac{\partial}{\partial t}-\Delta)\phi\leqslant \phi+C\phi^2-H.$$
With the description of $H$ before and the estimate for $\phi$, i.e., $\phi\leqslant C$, 
from the previous subsection, we can conclude that 
$$H\geqslant C|\nabla\phi|^2.$$ 
Now one arrives at $(\frac{\partial}{\partial t}-\Delta)\phi\leqslant C-C|\nabla\phi|^2$, 
which is also, for small enough positive constant $\epsilon$,  
$$(\frac{\partial}{\partial t}-\Delta)\phi +\epsilon |\nabla\phi|^2\leqslant C.$$

Of course, we also have
\begin{equation}
|(\nabla \phi, \frac{\nabla v}{C-v})|\leqslant \epsilon\cdot|\nabla \phi|^2+C\cdot|\nabla v|^2.
\end{equation}

Now consider the function $\Psi+\phi$. By choosing $\epsilon>0$ small enough above 
(which also affects the choices of $C$'s, but they will all be fixed eventually), we have 
$$(\frac{\partial}{\partial t}-\Delta)(\Psi+\phi)\leqslant C+C|\nabla v|^2-\epsilon\cdot|\nabla 
v|^4-(2-\epsilon){\rm Re}(\nabla(\Psi+\phi), \frac{\nabla v}{C-v}).$$

At the maximum value point of $\Psi+\phi$, we know $|\nabla v|^2$ can not be too large. 
It's then easy to conclude the upper bound for this term, and so for $\Psi$. Hence we have 
proved the gradient estimate, 
$$|\nabla v|\leqslant C.$$  

Now we want to do similar things for the Laplacian, $\Delta v$. Consider the quantity 
$\Phi=\frac{C-\Delta v}{C-v}$. Similar computation as before gives the following
\begin{equation}
\begin{split}
&~~ (\frac{\partial}{\partial t}-\Delta)\Phi\\
&= (\frac{\partial}{\partial t}-\Delta)(\frac{C-\Delta v}{C-v})\\
&= -\frac{1}{C-v}\cdot (\frac{\partial}{\partial t}-\Delta)\Delta v+\frac{C-\Delta v}{(C-v)^2}
\cdot(\frac{\partial}{\partial t}-\Delta)v\\
& ~~~~ +\frac{2{\rm Re}(\nabla v, \nabla\Delta v)}{(C-v)^2}-\frac{2|\nabla v|^2(C-\Delta 
v)}{(C-v)^3}\\
&= -\frac{1}{C-v}\cdot\bigl(\Delta v+({\rm Ric}(\tilde\omega_t), \sqrt{-1}\partial\bar\partial 
v)+\Delta \phi\bigr)+\frac{C-\Delta v}{C-v}\cdot(-n+\phi)\\
& ~~~~ +\frac{2{\rm Re}(\nabla v, \nabla\Delta v)}{(C-v)^2}-\frac{2|\nabla v|^2(C-\Delta 
v)}{(C-v)^3}. \nonumber
\end{split}
\end{equation}

We also have $\nabla(\frac{C-\Delta v}{C-v})=\frac{(C-\Delta v)\nabla v}{(C-v)^2}-
\frac{\nabla\Delta v}{C-v}$. Recall that it is already proved $(0\leqslant) \phi\leqslant 
C$. The following inequality follows from standard computation (as in \cite{song-tian}) 
and has actually been used before, 
$$\Delta \phi\geqslant ({\rm Ric}(\tilde\omega_t), \omega_T)+H-C\phi^2,$$
where $H\geqslant C|\nabla \phi|^2\geqslant 0$ from the bound of $\phi$ as mentioned 
before. Now we have 
$$({\rm Ric}, \sqrt{-1}\partial\bar\partial v)+\Delta \phi\geqslant ({\rm Ric}, \sqrt{-1}\partial
\bar\partial v+\omega_T)-C$$
with ${\rm Ric}={\rm Ric}(\tilde\omega_t)$. \\

At this moment, let's restrict to the case of $T=\infty$ as for Theorem 1.1. Then one has 
$$v=\frac{\partial u}{\partial t}+u, ~~~ {\rm Ric}=-\sqrt{-1}\partial\bar\partial v-
\omega_\infty.$$

Now the above estimate can be continued as follows,
\begin{equation}
\begin{split}
({\rm Ric}, \sqrt{-1}\partial\bar\partial v)+\Delta \phi
&\geqslant ({\rm Ric}, \sqrt{-1}\partial\bar\partial v+\omega_\infty)-C\\
&= -|\sqrt{-1}\partial\bar\partial v+\omega_\infty|^2-C\\
&\geqslant -(1+\epsilon)|\sqrt{-1}\partial\bar\partial v|^2-C\cdot|\omega_\infty|^2-C,
\nonumber
\end{split}
\end{equation}
where $\epsilon$ is a small positive number. As $\phi=\<\tilde\omega_t, \omega_\infty
\>\leqslant C$ and $\omega_\infty\geqslant 0$, which give $|\omega_\infty|^2\leqslant 
C$, one finally arrives at  
$$({\rm Ric}, \sqrt{-1}\partial\bar\partial v)+\Delta \phi\geqslant -(1+\epsilon)|\nabla\bar
\nabla v|^2-C.$$

As $T=\infty$, we also have $\Delta v=-R-\phi\leqslant C$, and so the numerator and 
denominator of $\Phi$ can both be positve and away from $0$. We only want to make 
sure that the numerator can not be too large. At this moment, we continue the 
computation for $\Phi$ as follows,
$$(\frac{\partial}{\partial t}-\Delta)\Phi\leqslant C+C\cdot(C-\Delta v)+\frac{(1+\epsilon)
|\nabla\bar\nabla v|^2}{C-v}-2{\rm Re}(\nabla \Phi, \frac{\nabla v}{C-v}).$$

Now we need the computation for $\Psi$ before, but from a different point of view. 
There is no need to involve $\epsilon$ there. Also, it is already known that $|\nabla v|
\leqslant C$. Basically from $(3.1)$ and $(3.2)$, we get the following which is useful 
here.
$$(\frac{\partial}{\partial t}-\Delta)(\Psi+\phi)\leqslant C-\frac{|\nabla\bar\nabla v|^2}
{C-v}-2{\rm Re}(\nabla(\Psi+\phi), \frac{\nabla v}{C-v}).$$
 
Combine them, by choosing $\epsilon$ small enough, to get 
$$(\frac{\partial}{\partial t}-\Delta)(\Phi+2\Psi+2\phi)\leqslant C+C\cdot(C-\Delta v)-
C\cdot |\nabla\bar\nabla v|^2-2{\rm Re}(\nabla(\Phi+2\Psi+2\phi), \frac{\nabla v}{C-v}).$$

The following observation would be enough to carry through the Maximum Principle 
argument,
$$|\nabla\bar\nabla v|^2\geqslant C(\Delta v)^2\geqslant C(C-\Delta v)^2-C.$$

Hence we can conclude that $\Phi\leqslant C$, and so
$$-\Delta v\leqslant C,$$
which is equivalent to $R\leqslant C$. Thus Theorem 1.1 is proved.

\begin{rema}

Of course, this would imply the bound of scalar curvature for the limiting metric (in 
the regular part), which is a rather trivial result as the metric is K\"ahler-Einstein 
there.  

In surface case, the limiting metric is an orbifold metric, which is more than just a 
metric with bounded scalar curvature over the regular part. For general dimension, 
hopefully this scalar curvature bound along K\"ahler-Ricci flow might be helpful to 
improve our knowledge on the limiting (singular K\"ahler-Einstein) metric.   

\end{rema}

\section{Remarks}

One might think of applying these computations for more general flows as introduced 
by H. Tsuji in \cite{tsu2} and discussed a little more in \cite{t-znote} and \cite{thesis}. 
More importantly, try to use them for this flow discussed above for finite time singularity 
case as the last part only works now for $T=\infty$. It seems to be quite interesting.

\end{document}